# Sharp Casorati Inequalities for Pointwise Bi-Slant Riemannian Maps into Sasakian Space Forms

Tanveer Fatima*

**Abstract**

We establish Casorati inequalities for pointwise bi-slant Riemannian maps into Sasakian space forms with the characteristic vector field tangent to the range. The Sasakian structure imposes specific conditions on the second fundamental form, leading to refined inequalities for the normalized horizontal scalar curvature that also account for O'Neill's tensor. For proper pointwise bi-slant maps, equality in the two principal inequalities is characterized by the vanishing of the second fundamental form on the orthogonal complement of the characteristic direction in the range. We also derive corresponding inequalities for pointwise semi-slant, hemi-slant, and slant Riemannian maps and construct examples attaining equality.

## 1 Introduction

Relations between intrinsic and extrinsic curvature invariants provide information about a submanifold that neither type of curvature alone can describe. One such extrinsic invariant is the Casorati curvature, defined in terms of the squared norm of the second fundamental form. Its normalized variants lead to inequalities involving scalar curvature and the Casorati curvatures of tangent hyperplanes. The equality cases carry further geometric information; see [1, 2, 3, 4].

The extension of these curvature inequalities to Riemannian maps requires a different geometric approach. Introduced by Fischer [5], Riemannian maps generalize both isometric immersions and Riemannian submersions by requiring the differential to preserve the metric on the orthogonal complement of its kernel. For a map of constant rank, the local factorization through its image allows us to relate the second fundamental form of the map to that of an immersed submanifold. For a broader discussion of Riemannian maps and their geometric properties, we refer to [6]. The horizontal scalar curvature of the source, however, is also affected by O'Neill's tensor $A$. Thus inequalities for Riemannian maps must distinguish the curvature of the image from the horizontal curvature of the source. Casorati inequalities for Riemannian maps have also been studied in [7, 8, 10, 11].

*Department of Mathematics and Statistics, College of Science in Yanbu, Taibah University, Yanbu Governorate, Kingdom of Saudi Arabia.
Email: `tansari@taibahu.edu.sa`

For Riemannian maps into Sasakian space forms, we refer to [10] and its subsequent corrigendum [12].

The geometry of slant distributions is characterized by the angle between the image of a tangent vector under the ambient structure tensor and the tangent space. The study of slant submanifolds [13] and slant Riemannian maps [14] has led to several related classes of maps. Pointwise slant Riemannian maps from Kähler manifolds were studied in [15], while Casorati inequalities along the vertical and horizontal distributions of pointwise slant Riemannian submersions were established in [9]. Pointwise semi-slant and hemi-slant Riemannian maps into almost Hermitian manifolds, together with associated curvature inequalities, were considered in [16, 17]. Pointwise bi-slant Riemannian maps from almost Hermitian manifolds were investigated in [18], whereas bi-slant Riemannian maps into Kenmotsu manifolds, including associated Casorati inequalities, were examined in [19]. Related estimates for bi-slant submanifolds and for Legendrian, slant, and bi-slant submanifolds of contact metric spaces appear in [20, 21, 22, 23, 24]. When the characteristic vector field is tangent to the range, the Sasakian structure determines the mixed components of the second fundamental form, which play an important role in the inequalities established below.

Let $F : (M_1, g_1) \to (M_2, \varphi_2, \xi_2, \eta_2, g_2)$ be a constant-rank Riemannian map into a Sasakian manifold, and suppose that the characteristic vector field $\xi_2$ is tangent to $\operatorname{Ran} F_*$. If $Q^{\mathcal{H}}U$ denotes the normal component of $\varphi_2 U$ and $h$ is the second fundamental form on horizontal vectors, the Sasakian identity $\nabla^{M_2}_U \xi_2 = -\varphi_2 U$ gives

$$h(U, \xi_2) = -Q^{\mathcal{H}}U, \qquad h(\xi_2, \xi_2) = 0.$$

For a pointwise bi-slant range decomposition

$$\operatorname{Ran} F_* = \mathcal{D}_1 \oplus \mathcal{D}_2 \oplus \langle \xi_2 \rangle,$$

put $\mathcal{D} = \mathcal{D}_1 \oplus \mathcal{D}_2$. If $\theta_1, \theta_2$ are the slant functions and $r_i = \dim \mathcal{D}_i$, then

$$\|h\|^2 = \|h|_{\mathcal{D}\times\mathcal{D}}\|^2 + 2r_1 \sin^2\theta_1 + 2r_2 \sin^2\theta_2.$$

For a proper pointwise bi-slant Riemannian map, the Sasakian identity determines nonzero mixed components of the second fundamental form, even when equality holds in the refined Casorati inequalities. These prescribed components are taken into account in deriving the estimates below.

To state our principal results, let $\mathcal{H} = (\ker F_*)^{\perp}$ and $r = \operatorname{rank} F$. We use the pointwise bi-slant decomposition above, with slant functions $\theta_i \in [0, \pi/2]$. The map is called *proper* when $0 < \theta_1, \theta_2 < \pi/2$. We write $h$ for the second fundamental form on horizontal vectors and $Q^{\mathcal{H}}$ for the normal component of $\varphi_2$ on $\mathcal{D}$. The quantities $\mathcal{C}^{\mathcal{H}}$, $\delta^{\mathcal{H}}_C(r-1)$, and $\widehat{\delta}^{\mathcal{H}}_C(r-1)$ denote the Casorati curvatures, while $\mathcal{C}(L)$ is the Casorati curvature of a hyperplane $L \subset \operatorname{Ran} F_*$. We denote by $\rho^{\mathcal{H}}$ and $\rho^{\mathcal{R}}$ the normalized horizontal and ambient curvature traces, respectively, and retain $A$ for O'Neill's tensor. Precise definitions are given in Sections 2 and 3.

We now state the two principal Casorati inequalities for pointwise bi-slant Riemannian maps into Sasakian space forms.

**Theorem A.** *Let $F : (M_1, g_1) \to (M_2(c), \varphi_2, \xi_2, \eta_2, g_2)$ be a pointwise bi-slant Riemannian map of constant rank $r \geq 3$ into a Sasakian space form, with $\xi_2 \in \Gamma(\operatorname{Ran} F_*)$.*

*Then*

$$\rho^{\mathcal{H}} \leq \rho^{\mathcal{R}} + \delta_C^{\mathcal{H}}(r-1) - \frac{\mathcal{C}^{\mathcal{H}}}{2(r-1)} - \frac{1}{2r}\inf_L \mathcal{C}(L)$$
$$- \frac{\|Q^{\mathcal{H}}\|^2}{r-1} - \frac{3\|A\|^2}{r(r-1)}. \tag{1.1}$$

*Moreover, under the condition* $\operatorname{rank}(Q^{\mathcal{H}}|_{\mathcal{D}}) \geq 2$*, equality holds if and only if*

$$h(U,V) = 0 \qquad \text{for all } U, V \in \mathcal{D}. \tag{1.2}$$

**Theorem B.** *Let* $F : (M_1, g_1) \to (M_2(c), \varphi_2, \xi_2, \eta_2, g_2)$ *be a proper pointwise bi-slant Riemannian map of constant rank* $r \geq 7$ *into a Sasakian space form, with* $\xi_2 \in \Gamma(\operatorname{Ran} F_*)$*. Then*

$$\rho^{\mathcal{H}} \leq \rho^{\mathcal{R}} + \widehat{\delta}_C^{\mathcal{H}}(r-1) - \frac{\mathcal{C}^{\mathcal{H}}}{2(r-1)} - \frac{2\|Q^{\mathcal{H}}\|^2}{r}$$
$$- \frac{2r-1}{r(r-1)}\min\{\sin^2\theta_1, \sin^2\theta_2\} - \frac{3\|A\|^2}{r(r-1)}. \tag{1.3}$$

*Equality holds at* $p \in M_1$ *if and only if*

$$h(U,V) = 0 \qquad \text{for all } U, V \in \mathcal{D}_p. \tag{1.4}$$

Both inequalities retain the contribution of O'Neill's tensor $A$ without requiring integrability of the horizontal distribution. Their proofs are given in Section 3.

For the second normalized Casorati curvature, a sharp estimate in rank five is established in Proposition 3.2. The corresponding inequalities for pointwise semi-slant, hemi-slant, and slant Riemannian maps are presented in three corollaries. Examples in Sasakian Heisenberg space forms demonstrate that the bounds are attained, including when the horizontal distribution is non-integrable.

The paper is organized as follows. Section 2 recalls the necessary background on Riemannian maps and Sasakian geometry and establishes the basic identities for pointwise bi-slant Riemannian maps. Section 3 introduces the Casorati curvatures and contains the proofs of Theorems A and B, the rank-five result, and the three corollaries. Section 4 provides examples attaining equality, and Section 5 concludes the paper.

# 2 Preliminaries

Throughout, all manifolds are assumed to be smooth and connected.

Let $F : (M_1^{m_1}, g_1) \to (M_2^{m_2}, g_2)$ be a smooth map between Riemannian manifolds. For a point $p \in M_1$, denote by $F_{*p} : T_pM_1 \to T_{F(p)}M_2$ the derivative map. Let $\ker F_{*p}$ be the kernel of $F_{*p}$ and $(\ker F_{*p})^\perp$ its orthogonal complement in $T_pM_1$. Then we have the orthogonal decomposition

$$T_pM_1 = \ker F_{*p} \oplus (\ker F_{*p})^\perp. \tag{2.1}$$

Similarly, at $F(p) \in M_2$, let $\operatorname{Ran} F_{*p}$ be the range of $F_{*p}$ and $(\operatorname{Ran} F_{*p})^\perp$ its orthogonal complement in $T_{F(p)}M_2$, giving

$$T_{F(p)}M_2 = \operatorname{Ran} F_{*p} \oplus (\operatorname{Ran} F_{*p})^\perp. \tag{2.2}$$

**Definition 2.1** ([5]). *$F$ is called a Riemannian map at $p \in M_1$ if the linear map $F_{*p} : (\ker F_{*p})^\perp \to \operatorname{Ran} F_{*p}$ is a linear isometry; that is,*

$$g_1(X, Y) = g_2(F_*X, F_*Y) \quad \text{for all } X, Y \in \Gamma((\ker F_*)^\perp). \tag{2.3}$$

If (2.3) holds at every point $p \in M_1$, then $F$ is called a Riemannian map.

When $(\operatorname{Ran} F_*)^\perp = \{0\}$, the map $F$ is a Riemannian submersion; when $\ker F_* = \{0\}$, it is an isometric immersion. Thus Riemannian maps simultaneously generalize isometric immersions and Riemannian submersions.

Let $\nabla^{M_1}$ and $\nabla^{M_2}$ denote the Levi-Civita connections on $M_1$ and $M_2$, respectively. The second fundamental form of $F$ is defined by [26]

$$(\nabla F_*)(X, Y) = \nabla^{M_2}_{F_*X} F_*Y - F_*(\nabla^{M_1}_X Y), \quad X, Y \in \Gamma(TM_1). \tag{2.4}$$

**Proposition 2.2** ([27]). *For a Riemannian map $F$, the second fundamental form satisfies*

$$(\nabla F_*)(X, Y) \in \Gamma((\operatorname{Ran} F_*)^\perp) \quad \textit{for all } X, Y \in \Gamma((\ker F_*)^\perp). \tag{2.5}$$

For $V \in \Gamma((\operatorname{Ran} F_*)^\perp)$, the shape operator $\mathcal{S}_V$ is defined by

$$\nabla^{M_2}_{F_*X} V = -\mathcal{S}_V F_*X + \nabla^{F\perp}_X V, \tag{2.6}$$

where $\nabla^{F\perp}$ is the connection in the normal bundle. The shape operator relates to the second fundamental form by

$$g_2((\nabla F_*)(X, Y), V) = g_2(\mathcal{S}_V F_*X, F_*Y). \tag{2.7}$$

**Proposition 2.3.** *Let $F : (M_1, g_1) \to (M_2, g_2)$ be a constant-rank Riemannian map. For horizontal vector fields $X, Y, Z, W$,*

$$\begin{aligned} g_2(R^{M_2}(F_*X, F_*Y)F_*Z, F_*W) &= g_1(R^{M_1}(X, Y)Z, W) \\ &+ g_2((\nabla F_*)(X, Z), (\nabla F_*)(Y, W)) - g_2((\nabla F_*)(X, W), (\nabla F_*)(Y, Z)) \\ &- 2g_1(A_X Y, A_Z W) - g_1(A_X Z, A_Y W) + g_1(A_Y Z, A_X W), \end{aligned} \tag{2.8}$$

*where $A_X Y = \mathcal{V}(\nabla^{M_1}_X Y)$ and $\mathcal{V}$ is the orthogonal projection onto $\ker F_*$. In particular, $A_X Y = -A_Y X$ for horizontal fields.*

*Proof.* Locally, $F = i \circ \pi$, where $\pi : M_1 \to S$ is a Riemannian submersion onto the image equipped with its induced metric, and $i : S \to M_2$ is an isometric immersion. For horizontal fields, O'Neill's curvature identity [28] gives

$$\begin{aligned} g_S(R^S(\pi_*X, \pi_*Y)\pi_*Z, \pi_*W) &= g_1(R^{M_1}(X, Y)Z, W) \\ &- 2g_1(A_X Y, A_Z W) - g_1(A_X Z, A_Y W) + g_1(A_Y Z, A_X W). \end{aligned}$$

The Gauss equation for $i$, together with the identification of its second fundamental form with $(\nabla F_*)$ on horizontal vectors, yields (2.8). Moreover, for horizontal fields, $A_X Y = \frac{1}{2}\mathcal{V}[X, Y]$; hence $A$ is skew-symmetric. □

*Remark* 2.4. The horizontal curvature identity used in [25] does not display the $A$-tensor contribution. In its traced form, this contribution is $-3\|A\|^2$. Consequently, the upper Casorati bounds stated there remain valid, but attainment of equality for a general Riemannian map requires $A = 0$ in addition to the algebraic equality conditions for these curvature inequalities.

We shall use the following elementary norm identity. If $\{e_1, \ldots, e_r\}$ is an orthonormal basis of an inner product space $V$ and $T : V \to V$ is a linear operator, then

$$\|T\|^2 := \sum_{i,j=1}^{r} g(Te_i, e_j)^2 = \sum_{j=1}^{r} \|Te_j\|^2.$$

This is the usual Parseval identity.

Recall that an almost contact metric manifold $(M_2, \varphi_2, \xi_2, \eta_2, g_2)$ satisfies

$$\varphi_2^2 U = -U + \eta_2(U)\xi_2, \qquad \eta_2(\xi_2) = 1, \qquad \eta_2(U) = g_2(U, \xi_2),$$
$$g_2(\varphi_2 U, \varphi_2 V) = g_2(U, V) - \eta_2(U)\eta_2(V).$$

The manifold is Sasakian if

$$(\nabla_U^{M_2} \varphi_2)V = g_2(U, V)\xi_2 - \eta_2(V)U.$$

In particular, $\nabla_U^{M_2} \xi_2 = -\varphi_2 U$. We use the curvature convention $R(U, V)W = \nabla_U \nabla_V W - \nabla_V \nabla_U W - \nabla_{[U,V]} W$. If $M_2 = M_2(c)$ has constant $\varphi_2$-sectional curvature $c$, its curvature tensor is

$$\begin{aligned} R^{M_2}(U, V)W = {} & \frac{c+3}{4}\{g_2(V, W)U - g_2(U, W)V\} \\ & + \frac{c-1}{4}\{g_2(U, \varphi_2 W)\varphi_2 V - g_2(V, \varphi_2 W)\varphi_2 U \\ & \qquad + 2g_2(U, \varphi_2 V)\varphi_2 W\} \\ & + \frac{c-1}{4}\{\eta_2(U)\eta_2(W)V - \eta_2(V)\eta_2(W)U \\ & \qquad + g_2(U, W)\eta_2(V)\xi_2 - g_2(V, W)\eta_2(U)\xi_2\}. \end{aligned} \tag{2.9}$$

See, for example, [29, 30].

We consider pointwise bi-slant Riemannian maps into almost contact metric manifolds with the characteristic vector field tangent to the range. Throughout, $F$ has constant rank $r$. All decompositions are understood along $F$.

Let $F : (M_1, g_1) \to (M_2, \varphi_2, \xi_2, \eta_2, g_2)$ be a Riemannian map into an almost contact metric manifold, and assume that $\xi_2 \in \Gamma(\operatorname{Ran} F_*)$. Put

$$\mathcal{D} = \operatorname{Ran} F_* \cap \ker \eta_2.$$

For $U \in \Gamma(\mathcal{D})$, we write the tangential and normal components of $\varphi_2 U$, following [25], as

$$\varphi_2 U = P^{\mathcal{H}} U + Q^{\mathcal{H}} U, \tag{2.10}$$

where $P^{\mathcal{H}} U \in \Gamma(\operatorname{Ran} F_*)$ and $Q^{\mathcal{H}} U \in \Gamma((\operatorname{Ran} F_*)^\perp)$. The superscript $\mathcal{H}$ reflects the identification of the horizontal distribution with $\operatorname{Ran} F_*$ through $F_*$. Since $\varphi_2 \xi_2 = 0$, we have $g_2(P^{\mathcal{H}} U, \xi_2) = 0$ for every $U \in \mathcal{D}$. Hence $P^{\mathcal{H}}(\mathcal{D}) \subseteq \mathcal{D}$.

For a local orthonormal frame $\{E_1, \ldots, E_{r-1}\}$ of $\mathcal{D}$, we write

$$\|P^{\mathcal{H}}\|^2 = \sum_{i=1}^{r-1} \|P^{\mathcal{H}} E_i\|^2, \qquad \|Q^{\mathcal{H}}\|^2 = \sum_{i=1}^{r-1} \|Q^{\mathcal{H}} E_i\|^2.$$

**Definition 2.5.** A Riemannian map $F : (M_1, g_1) \to (M_2, \varphi_2, \xi_2, \eta_2, g_2)$, with $\xi_2 \in \Gamma(\operatorname{Ran} F_*)$, is called a *pointwise bi-slant Riemannian map* if

$$\operatorname{Ran} F_* = \mathcal{D}_1 \oplus \mathcal{D}_2 \oplus \langle \xi_2 \rangle \tag{2.11}$$

is an orthogonal decomposition into nonzero distributions $\mathcal{D}_1, \mathcal{D}_2 \subseteq \ker \eta_2$ and the Reeb direction, such that $\varphi_2 \mathcal{D}_1 \perp \mathcal{D}_2$ and there exist functions $\theta_1, \theta_2 : M_1 \to [0, \pi/2]$ satisfying

$$\cos \theta_i(p) = \frac{\|P^{\mathcal{H}} U\|}{\|U\|}, \qquad 0 \neq U \in (\mathcal{D}_i)_p, \quad i = 1, 2. \tag{2.12}$$

The decomposition (2.11) determines the distribution

$$\mathcal{D} = \mathcal{D}_1 \oplus \mathcal{D}_2 = \operatorname{Ran} F_* \cap \ker \eta_2,$$

which is orthogonal to $\xi_2$. The map $F$ is called *proper pointwise bi-slant* if $0 < \theta_1, \theta_2 < \pi/2$.

We write $r_i = \dim \mathcal{D}_i$ ($i = 1, 2$). On $\mathcal{D}$, the almost contact metric identities give $\varphi_2^2 = -I$ and $g_2(\varphi_2 U, \varphi_2 V) = g_2(U, V)$.

**Proposition 2.6.** *Let $F$ be a pointwise bi-slant Riemannian map into an almost contact metric manifold with $\xi_2 \in \Gamma(\operatorname{Ran} F_*)$. For $U, V \in \Gamma(\mathcal{D}_i)$, $i = 1, 2$, we have*

$$(P^{\mathcal{H}})^2 U = -\cos^2 \theta_i \, U, \tag{2.13}$$

$$g_2(P^{\mathcal{H}} U, P^{\mathcal{H}} V) = \cos^2 \theta_i \, g_2(U, V), \tag{2.14}$$

$$g_2(Q^{\mathcal{H}} U, Q^{\mathcal{H}} V) = \sin^2 \theta_i \, g_2(U, V). \tag{2.15}$$

*Consequently,*

$$\|P^{\mathcal{H}}\|^2 = r_1 \cos^2 \theta_1 + r_2 \cos^2 \theta_2, \tag{2.16}$$

$$\|Q^{\mathcal{H}}\|^2 = r_1 \sin^2 \theta_1 + r_2 \sin^2 \theta_2. \tag{2.17}$$

*Proof.* The orthogonality condition implies that $P^{\mathcal{H}}$ preserves each $\mathcal{D}_i$. The slant identity and polarization give $g_2(P^{\mathcal{H}} U, P^{\mathcal{H}} V) = \cos^2 \theta_i \, g_2(U, V)$. Since $P^{\mathcal{H}}$ is skew-symmetric, this is equivalent to (2.13). Using $g_2(\varphi_2 U, \varphi_2 V) = g_2(U, V)$ gives (2.15). Taking traces and using $\varphi_2 \xi_2 = 0$ yields the norm identities. □

*Remark* 2.7. For a proper pointwise bi-slant map, both $r_1$ and $r_2$ are even, since $P^{\mathcal{H}} / \cos \theta_i$ defines an almost complex structure on $\mathcal{D}_i$. Consequently,

$$r = r_1 + r_2 + 1 \geq 5, \qquad r \text{ is odd.} \tag{2.18}$$

When the target is Sasakian, the Reeb vector field also restricts the second fundamental form in a useful way. Since $F$ has constant rank, its image is locally an immersed submanifold $S$ of $M_2$. We write $h(U, V)$ for the second fundamental form of $S$, identified with $(\nabla F_*)$ on horizontal lifts of $U, V \in \Gamma(\operatorname{Ran} F_*)$.

Let $r = \operatorname{rank} F$ and $B = h|_{\mathcal{D} \times \mathcal{D}}$. Choose local orthonormal frames $\{e_1, \ldots, e_r\}$ of $\operatorname{Ran} F_*$, with $e_r = \xi_2$, and $\{E_1, \ldots, E_{r-1}\}$ of $\mathcal{D}$. We use the conventions

$$\|h\|^2 = \sum_{a,b=1}^{r} \|h(e_a, e_b)\|^2, \qquad \operatorname{tr} h = \sum_{a=1}^{r} h(e_a, e_a),$$

$$\|B\|^2 = \sum_{i,j=1}^{r-1} \|B(E_i, E_j)\|^2, \quad \operatorname{tr} B = \sum_{i=1}^{r-1} B(E_i, E_i).$$

**Proposition 2.8.** *Let $F : (M_1, g_1) \to (M_2, \varphi_2, \xi_2, \eta_2, g_2)$ be a constant-rank Riemannian map into a Sasakian manifold, with $\xi_2 \in \Gamma(\operatorname{Ran} F_*)$.*

*Then, for $U \in \Gamma(\mathcal{D})$,*

$$h(U, \xi_2) = -Q^{\mathcal{H}} U, \qquad h(\xi_2, \xi_2) = 0. \tag{2.19}$$

*Moreover,*

$$\|h\|^2 = \|B\|^2 + 2\|Q^{\mathcal{H}}\|^2, \tag{2.20}$$

$$\operatorname{tr} h = \operatorname{tr} B. \tag{2.21}$$

*In particular, if $F$ is pointwise bi-slant, then*

$$\|h\|^2 = \|B\|^2 + 2r_1 \sin^2\theta_1 + 2r_2 \sin^2\theta_2. \tag{2.22}$$

*Proof.* Taking the normal component of the Sasakian identity $\nabla^{M_2}_U \xi_2 = -\varphi_2 U$ gives $h(U, \xi_2) = -Q^{\mathcal{H}} U$. Since $\varphi_2 \xi_2 = 0$, we also have $h(\xi_2, \xi_2) = 0$.

Choose an orthonormal frame $\{E_1, \dots, E_{r-1}\}$ of $\mathcal{D}$. Together with $\xi_2$, it forms an orthonormal frame of $\operatorname{Ran} F_*$. Consequently,

$$\|h\|^2 = \|B\|^2 + 2\sum_{i=1}^{r-1} \|Q^{\mathcal{H}} E_i\|^2 = \|B\|^2 + 2\|Q^{\mathcal{H}}\|^2.$$

The trace identity follows from $h(\xi_2, \xi_2) = 0$. For a pointwise bi-slant map, (2.17) gives (2.22). □

*Remark* 2.9. The assumption $\xi_2 \in \Gamma(\operatorname{Ran} F_*)$ is essential in the results below. If $\xi_2$ is normal to the local image in a Sasakian manifold, the tangential part of $\nabla^{M_2}_U \xi_2 = -\varphi_2 U$ gives $\mathcal{S}_{\xi_2} U = P^{\mathcal{H}} U$. The shape operator is symmetric, whereas $P^{\mathcal{H}}$ is skew-symmetric. Thus $P^{\mathcal{H}} = 0$, and the image is anti-invariant. The case in which $\xi_2$ has both tangential and normal components is not considered here.

*Remark* 2.10. For a proper pointwise bi-slant Riemannian map, $Q^{\mathcal{H}} U \neq 0$ for every nonzero $U \in \mathcal{D}$. By (2.19), the local image therefore cannot be totally geodesic in the Sasakian target.

The two-angle formulation includes the familiar mixed classes as limiting cases. When $\theta_1 = 0$, the first summand is invariant and the map is pointwise semi-slant; when $\theta_1 = \pi/2$, it is anti-invariant and the map is pointwise hemi-slant. If $\theta_1 = \theta_2 = \theta$, the range is pointwise slant, as in [25]. Pointwise semi-slant and hemi-slant Riemannian maps into almost Hermitian manifolds were studied in [16, 17]. These special cases will be considered after the main theorems.

Finally, we record the ambient curvature terms which will enter the Casorati inequalities. Let $\{e_1, \dots, e_r\}$ be a local orthonormal frame of $\operatorname{Ran} F_*$, and define

$$\tau^{\mathcal{R}} = \sum_{1 \le a < b \le r} g_2\Big(R^{M_2}(e_a, e_b)e_b, e_a\Big). \tag{2.23}$$

This is the ambient sectional-curvature trace over the range distribution, rather than the intrinsic scalar curvature of the immersed image.

**Proposition 2.11.** *Let $F$ be a pointwise bi-slant Riemannian map of constant rank $r$ into a Sasakian space form $M_2(c)$, with $\xi_2 \in \Gamma(\operatorname{Ran} F_*)$. Then*

$$\begin{aligned} 2\tau^{\mathcal{R}} = & \frac{r(r-1)(c+3)}{4} + \frac{3(c-1)}{4}\Big(r_1 \cos^2\theta_1 + r_2 \cos^2\theta_2\Big) \\ & - \frac{(r-1)(c-1)}{2}. \end{aligned} \tag{2.24}$$

*Proof.* Tracing (2.9) over an orthonormal frame of $\operatorname{Ran} F_*$, with $\xi_2$ as its last member, gives

$$2\tau^{\mathcal{R}} = \frac{r(r-1)(c+3)}{4} + \frac{3(c-1)}{4}\|P^{\mathcal{H}}\|^2 - \frac{(r-1)(c-1)}{2}.$$

Now apply (2.16). □

# 3 Casorati Inequalities for Pointwise Bi-Slant Riemannian Maps

Let $F : (M_1, g_1) \to (M_2, \varphi_2, \xi_2, \eta_2, g_2)$ be a constant-rank Riemannian map into a Sasakian manifold, with $r = \operatorname{rank} F \geq 3$ and $\xi_2 \in \Gamma(\operatorname{Ran} F_*)$. We identify the horizontal distribution with $\operatorname{Ran} F_*$ through $F_*$ and use the notation $\mathcal{D} = \operatorname{Ran} F_* \cap \ker \eta_2$ and $B = h|_{\mathcal{D}\times\mathcal{D}}$ introduced in Section 2.

For a hyperplane $L \subset \operatorname{Ran} F_*$ with orthonormal basis $\{\ell_1, \ldots, \ell_{r-1}\}$, we write

$$\|h|_{L\times L}\|^2 = \sum_{i,j=1}^{r-1} \|h(\ell_i, \ell_j)\|^2.$$

The Casorati curvature of the horizontal distribution and its restriction to a hyperplane $L \subset \operatorname{Ran} F_*$ are

$$\mathcal{C}^{\mathcal{H}} = \frac{\|h\|^2}{r}, \qquad \mathcal{C}(L) = \frac{\|h|_{L\times L}\|^2}{r-1}. \tag{3.1}$$

The corresponding normalized Casorati curvatures are

$$\delta_C^{\mathcal{H}}(r-1) = \frac{1}{2}\mathcal{C}^{\mathcal{H}} + \frac{r+1}{2r}\inf_L \mathcal{C}(L), \tag{3.2}$$

$$\widehat{\delta}_C^{\mathcal{H}}(r-1) = 2\mathcal{C}^{\mathcal{H}} - \frac{2r-1}{2r}\sup_L \mathcal{C}(L). \tag{3.3}$$

The extrema are taken over hyperplanes of $\operatorname{Ran} F_*$ and are attained since the Grassmannian of hyperplanes is compact.

Let $\{X_1, \ldots, X_r\}$ be an orthonormal horizontal frame and $e_i = F_* X_i$ the corresponding orthonormal frame of $\operatorname{Ran} F_*$. The horizontal scalar curvature is given by

$$\tau^{\mathcal{H}} = \sum_{1\leq i<j\leq r} g_1(R^{M_1}(X_i, X_j)X_j, X_i).$$

Let $\tau^{\mathcal{R}}$ be the ambient curvature trace of $M_2$ over $\operatorname{Ran} F_*$ defined in (2.23). Put

$$\rho^{\mathcal{H}} = \frac{2\tau^{\mathcal{H}}}{r(r-1)}, \qquad \rho^{\mathcal{R}} = \frac{2\tau^{\mathcal{R}}}{r(r-1)}.$$

For a pointwise bi-slant Riemannian map into a Sasakian space form $M_2(c)$, the ambient trace is

$$\rho^{\mathcal{R}} = \frac{c+3}{4} + \frac{3(c-1)}{4r(r-1)}\Big(r_1\cos^2\theta_1 + r_2\cos^2\theta_2\Big) - \frac{c-1}{2r}, \tag{3.4}$$

as follows from Proposition 2.11. For the same frame, write $\|A\|^2 = \sum_{i,j=1}^{r}\|A_{X_i}X_j\|^2$. Tracing the Gauss equation (2.8), we obtain

$$r(r-1)\rho^{\mathcal{H}} = r(r-1)\rho^{\mathcal{R}} + \|\operatorname{tr} h\|^2 - \|h\|^2 - 3\|A\|^2. \tag{3.5}$$

*Proof of Theorem A.* Fix a point and a hyperplane $L \subset \operatorname{Ran} F_*$. Choose $E_1, \ldots, E_{r-1}$ orthonormal in $\mathcal{D}$, with $E_1, \ldots, E_{r-2} \in L$. Set $v_i = h(E_i, E_i)$. From (2.21) and weighted Cauchy–Schwarz,

$$\begin{aligned} \|\operatorname{tr} h\|^2 = \left\|\sum_{i=1}^{r-1} v_i\right\|^2 &\le \frac{r}{2}\left(2\sum_{i=1}^{r-2}\|v_i\|^2 + \|v_{r-1}\|^2\right) \\ &\le \frac{r}{2}\left(\left\|h|_{\mathcal{D}\times\mathcal{D}}\right\|^2 + \|h|_{L\times L}\|^2\right). \end{aligned} \tag{3.6}$$

Combining (3.6) with (3.5) and (2.22) yields

$$\begin{aligned} \rho^{\mathcal{H}} \le \rho^{\mathcal{R}} + \frac{r-2}{2(r-1)}\mathcal{C}^{\mathcal{H}} + \frac{1}{2}\mathcal{C}(L) \\ - \frac{\|Q^{\mathcal{H}}\|^2}{r-1} - \frac{3\|A\|^2}{r(r-1)}. \end{aligned}$$

Since $L$ is arbitrary, we take the infimum and use

$$\frac{r-2}{2(r-1)}\mathcal{C}^{\mathcal{H}} + \frac{1}{2}\inf_L \mathcal{C}(L) = \delta_C^{\mathcal{H}}(r-1) - \frac{\mathcal{C}^{\mathcal{H}}}{2(r-1)} - \frac{1}{2r}\inf_L \mathcal{C}(L).$$

Suppose $\operatorname{rank}(Q^{\mathcal{H}}|_{\mathcal{D}}) \ge 2$ and equality holds. Take a hyperplane $L$ attaining the infimum, and write $L = \operatorname{span}\{E_1, \ldots, E_{r-2}, W\}$ with $W = \gamma E_{r-1} + \beta\xi_2$ and $\gamma^2 + \beta^2 = 1$. Equality in the second inequality of (3.6) forces $h(E_i, E_j) = 0$ for $i \ne j$, $1 \le i, j \le r-1$, and $h(E_i, W) = 0$ for $1 \le i \le r-2$. Thus

$$0 = h(E_i, W) = -\beta Q^{\mathcal{H}} E_i.$$

Since $E_1, \ldots, E_{r-2}$ span a hyperplane of $\mathcal{D}$ and $\operatorname{rank}(Q^{\mathcal{H}}|_{\mathcal{D}}) \ge 2$, we have $Q^{\mathcal{H}}E_i \ne 0$ for some $i \le r-2$. Hence $\beta = 0$ and $L = \mathcal{D}$. Equality in (3.6) then gives $v_{r-1} = 0$, while equality in weighted Cauchy–Schwarz implies $v_1 = \cdots = v_{r-2} = 0$. Hence (1.2) holds. Conversely, if (1.2) holds, take $L = \mathcal{D}$; equality holds in (3.6) and therefore in (1.1). □

For the second normalized Casorati curvature, the supremum over hyperplanes leads to a different optimization. In the proper pointwise bi-slant case, the nonvanishing normal components $Q^{\mathcal{H}}U$ yield an additional correction term.

**Lemma 3.1.** *Let $F : (M_1, g_1) \to (M_2, \varphi_2, \xi_2, \eta_2, g_2)$ be a Riemannian map of constant rank $r \ge 3$ into a Sasakian manifold, with $\xi_2 \in \Gamma(\operatorname{Ran} F_*)$. Put $B = h|_{\mathcal{D}\times\mathcal{D}}$, and*

$$B^\circ(U, V) = B(U, V) - \frac{g_2(U, V)}{r-1}\operatorname{tr} B.$$

*For any hyperplane $L \subset \operatorname{Ran} F_*$, write a unit normal to $L$ as $\sqrt{t}\,\xi_2 + \sqrt{1-t}\,u$, where $u \in \mathcal{D}$ is a unit vector and $0 \le t \le 1$, and extend $E_1 = u$ to an orthonormal frame $E_1, \dots, E_{r-1}$ of $\mathcal{D}$. Then*

$$
\begin{aligned}
\|h\|^2 - \|h|_{L\times L}\|^2 = {} & \left\|(1-t)B(u,u) - 2\sqrt{t(1-t)}Q^{\mathcal{H}}u\right\|^2 \\
& + 2\left\|\sqrt{t(1-t)}B(u,u) + (1-2t)Q^{\mathcal{H}}u\right\|^2 \\
& + 2\sum_{j=2}^{r-1}\left\|\sqrt{1-t}B(u,E_j) - \sqrt{t}\,Q^{\mathcal{H}}E_j\right\|^2 .
\end{aligned}
\tag{3.7}
$$

*Moreover,*

$$
\begin{aligned}
& \frac{2r-1}{2}\left(\|h\|^2 - \|h|_{L\times L}\|^2\right) + (r-1)\|B^\circ\|^2 \\
& \quad \ge (2r-1)\left\{\frac{1-t}{1+t}\|Q^{\mathcal{H}}u\|^2 + \frac{(r-1)t}{(r-1)+\frac{2r-1}{2}(1-t)}\left(\|Q^{\mathcal{H}}\|^2 - \|Q^{\mathcal{H}}u\|^2\right)\right\}.
\end{aligned}
\tag{3.8}
$$

*For $0 \le t < 1$, equality in the last estimate holds exactly when, with $q_j = Q^{\mathcal{H}}E_j$ and $k = (2r-1)/2$,*

$$
\begin{aligned}
B(E_i,E_i) &= \frac{2t^{3/2}}{(1+t)\sqrt{1-t}}q_1, \qquad 1 \le i \le r-1, \\
B(u,E_j) &= \frac{k\sqrt{t(1-t)}}{(r-1)+k(1-t)}q_j, \quad 2 \le j \le r-1,
\end{aligned}
$$

*and all remaining off-diagonal components of $B$ vanish.*

*Proof.* The unit vector $\sqrt{t}\,u - \sqrt{1-t}\,\xi_2$ lies in $L$. The identity (3.7) follows by expanding $h$ on the resulting orthonormal frame and using (2.19). It remains to establish (3.8).

Put $k = (2r-1)/2$ and $q_j = Q^{\mathcal{H}}E_j$. The trace-free identity

$$
(r-1)\|B^\circ\|^2 = \sum_{i<j}\|B(E_i,E_i) - B(E_j,E_j)\|^2 + 2(r-1)\sum_{i<j}\|B(E_i,E_j)\|^2.
$$

Since (3.7) depends on the diagonal components of $B$ only through $B(u,u)$, the non-negative sum of squared differences is minimized, for fixed $B(u,u)$, when all diagonal components coincide. Similarly, the off-diagonal components $B(E_i,E_j)$, $2 \le i < j \le r-1$, may be taken to be zero.

Write their common diagonal value as $\Lambda$ and put $Y_j = B(u,E_j)$, $2 \le j \le r-1$. The remaining quadratic expression is

$$
\begin{aligned}
& k\Big\{(1-t)(1+t)\|\Lambda\|^2 - 4t\sqrt{t(1-t)}\langle\Lambda,q_1\rangle + 2(1-2t+2t^2)\|q_1\|^2\Big\} \\
& \quad + 2\sum_{j=2}^{r-1}\Big\{[(r-1)+k(1-t)]\|Y_j\|^2 - 2k\sqrt{t(1-t)}\langle Y_j,q_j\rangle + kt\|q_j\|^2\Big\}.
\end{aligned}
$$

For $0 \le t < 1$, completion of squares gives the minimizing values

$$
\Lambda = \frac{2t^{3/2}}{(1+t)\sqrt{1-t}}q_1, \qquad Y_j = \frac{k\sqrt{t(1-t)}}{(r-1)+k(1-t)}q_j.
$$

Substitution yields

$$2k\frac{1-t}{1+t}\|q_1\|^2 + 2k\frac{(r-1)t}{(r-1)+k(1-t)}\sum_{j=2}^{r-1}\|q_j\|^2,$$

which is precisely the right-hand side of (3.8). The inequality at $t = 1$ follows by continuity. The equality conditions follow from the same completed squares. $\square$

*Proof of Theorem B.* Write $B = h|_{\mathcal{D}\times\mathcal{D}}$ and let $B^\circ$ be as in Lemma 3.1. Let $L$ be a hyperplane at which $\mathcal{C}(L)$ is maximal. By (3.5), (3.3), and (2.22),

$$\begin{aligned}\rho^{\mathcal{H}} &= \rho^{\mathcal{R}} + \widehat{\delta}^{\mathcal{H}}_C(r-1) - \frac{\mathcal{C}^{\mathcal{H}}}{2(r-1)} - \frac{2\|Q^{\mathcal{H}}\|^2}{r} - \frac{3\|A\|^2}{r(r-1)} \\ &\quad - \frac{1}{r(r-1)}\left\{\frac{2r-1}{2}\left(\|h\|^2 - \|h|_{L\times L}\|^2\right) + (r-1)\|B^\circ\|^2\right\}. \qquad (3.9)\end{aligned}$$

Here we used $(r-1)\|B^\circ\|^2 = (r-1)\|B\|^2 - \|\operatorname{tr} B\|^2$. Lemma 3.1 now gives the common quadratic estimate (3.8). It remains to use the rank and slant assumptions.

Since

$$\|Q^{\mathcal{H}}u\|^2 \geq \min\{\sin^2\theta_1, \sin^2\theta_2\}, \qquad \|Q^{\mathcal{H}}\|^2 - \|Q^{\mathcal{H}}u\|^2 \geq (r-2)\min\{\sin^2\theta_1, \sin^2\theta_2\},$$

we need only observe that, for $r \geq 7$ and $0 \leq t \leq 1$,

$$\frac{1-t}{1+t} + \frac{(r-1)(r-2)t}{(r-1)+\frac{2r-1}{2}(1-t)} \geq 1.$$

For $t > 0$, this inequality is equivalent to

$$r^2 - 7r + 5 + (r^2 - r + 1)t \geq 0,$$

which holds for $r \geq 7$. At $t = 0$, the original inequality is an equality. Substitution in (3.9) proves (1.3).

For equality, the last estimate requires $t = 0$. Equality in the completed squares then forces $B = 0$. Conversely, if $B = 0$, take $u$ in a slant distribution corresponding to the smaller value of $\sin^2\theta_i$. The hyperplane normal to $u$ maximizes $\mathcal{C}(L)$: indeed, for any hyperplane with unit normal $\sqrt{t}\,\xi_2 + \sqrt{1-t}\,v$, (3.7) with $B = 0$ gives

$$\begin{aligned}\|h\|^2 - \|h|_{L\times L}\|^2 &= 2(1 - 2t + 2t^2)\|Q^{\mathcal{H}}v\|^2 + 2t(\|Q^{\mathcal{H}}\|^2 - \|Q^{\mathcal{H}}v\|^2) \\ &\geq 2\min\{\sin^2\theta_1, \sin^2\theta_2\}.\end{aligned}$$

Here the last inequality uses $r \geq 4$. The chosen hyperplane attains this lower bound, so equality holds in (3.8) and (3.9). This proves (1.4). $\square$

The smallest possible rank of a proper pointwise bi-slant map in the present setting is five. The same argument gives a more precise correction in this case, depending on the two slant functions.

**Proposition 3.2.** *Let $F$ satisfy the assumptions of Theorem B, except that $r = 5$. Relabel the distributions, if necessary, so that $\sin^2\theta_1 \leq \sin^2\theta_2$. If*

$$\sin^2\theta_2 \geq \frac{13}{8}\sin^2\theta_1,$$

*then* (1.3) *remains valid with* $r = 5$*, and equality holds if and only if* $h|_{\mathcal{D}\times\mathcal{D}} = 0$*. If*

$$\sin^2\theta_2 < \frac{13}{8}\sin^2\theta_1,$$

*then* (1.3) *holds with the term* $-\frac{9}{20}\sin^2\theta_1$ *replaced by*

$$-\frac{9}{65}\left[\sqrt{17\sin^2\theta_1\big(\sin^2\theta_1 + 2\sin^2\theta_2\big)} - 2\big(\sin^2\theta_1 + \sin^2\theta_2\big)\right]. \tag{3.10}$$

*In this second case, equality at a point holds precisely when, for the minimizing value*

$$t = \frac{17 - \sqrt{68(1 + 2\sin^2\theta_2/\sin^2\theta_1)}}{9 + \sqrt{68(1 + 2\sin^2\theta_2/\sin^2\theta_1)}},$$

*there is a unit vector* $u \in \mathcal{D}_1$ *if* $\sin^2\theta_1 < \sin^2\theta_2$*, or an arbitrary unit vector* $u \in \mathcal{D}$ *if the angles agree, and an orthonormal frame* $E_1 = u, E_2, E_3, E_4$ *of* $\mathcal{D}$*, for which*

$$B(E_i, E_i) = \frac{2t^{3/2}}{(1+t)\sqrt{1-t}}Q^{\mathcal{H}}u, \qquad 1 \le i \le 4, \tag{3.11}$$

$$B(u, E_j) = \frac{9\sqrt{t(1-t)}}{17 - 9t}Q^{\mathcal{H}}E_j, \qquad 2 \le j \le 4, \tag{3.12}$$

*with all other off-diagonal components zero, where* $B = h|_{\mathcal{D}\times\mathcal{D}}$*.*

*Proof.* For $r = 5$, the proper slant distributions both have rank two. In (3.8), write temporarily $\mu = \sin^2\theta_1$, $\nu = \sin^2\theta_2$, and $\lambda = \|Q^{\mathcal{H}}u\|^2$, so that $\mu \le \lambda \le \nu$ and $\|Q^{\mathcal{H}}\|^2 = 2(\mu+\nu)$. The right-hand side of (3.8) is

$$9\left\{\frac{1-t}{1+t}\lambda + \frac{8t}{17-9t}\big(2(\mu+\nu) - \lambda\big)\right\}. \tag{3.13}$$

For $t \le 17 - 4\sqrt{17}$ the coefficient of $\lambda$ is nonnegative, so (3.13) is minimized at $\lambda = \mu$. For larger $t$, the minimizing value is $\lambda = \nu$. In this range, however, $\frac{1-t}{1+t} + \frac{24t}{17-9t} \ge 1$; thus the expression is at least $9\mu$. If $\mu = \nu$, every unit $u \in \mathcal{D}$ gives $\lambda = \mu$; otherwise $\lambda = \mu$ requires $u \in \mathcal{D}_1$. It remains to minimize

$$9\left\{\frac{1-t}{1+t}\mu + \frac{8t}{17-9t}(\mu + 2\nu)\right\}.$$

This function is strictly convex on $0 \le t \le 1$, and its derivative at $t = 0$ is nonnegative exactly when $\nu \ge 13\mu/8$. In that case the minimum is $9\mu$, attained at $t = 0$. Otherwise the minimum occurs at

$$t = \frac{17 - \sqrt{68(1 + 2\nu/\mu)}}{9 + \sqrt{68(1 + 2\nu/\mu)}}$$

and is

$$\frac{36}{13}\left[\sqrt{17\mu(\mu + 2\nu)} - 2(\mu+\nu)\right].$$

Substituting these values into (3.9) gives the stated inequalities. In the first case, equality forces $t = 0$ and $B = 0$, and the converse follows as in Theorem B.

For equality in the second case, the quadratic minimizers above give (3.11) and (3.12). Conversely, these conditions make (3.8) an equality for the hyperplane normal to $\sqrt{t}\,\xi_2 + \sqrt{1-t}\,u$. As the lower bound holds for every hyperplane, this one maximizes $\mathcal{C}(L)$, and equality follows from (3.9). □

*Remark* 3.3. The estimate for $\widehat{\delta}^{\mathcal{H}}_C(r-1)$ in Theorem B requires $r \geq 7$; the admissible rank-five case is treated separately in Proposition 3.2. For hemi-slant maps, rank six is not covered by this optimization argument.

The two principal estimates give three specializations. Each corollary states the $\delta^{\mathcal{H}}_C$ inequality and, under its additional rank and slant assumptions, the $\widehat{\delta}^{\mathcal{H}}_C$ inequality. Throughout, $\rho^{\mathcal{R}}$ is the ambient trace in (3.4), with the corresponding specialization of the slant functions.

**Corollary 3.4.** *Let $F : (M_1, g_1) \to (M_2(c), \varphi_2, \xi_2, \eta_2, g_2)$ be a pointwise semi-slant Riemannian map of constant rank $r \geq 3$ into a Sasakian space form, with $\xi_2 \in \Gamma(\operatorname{Ran} F_*)$ and $\operatorname{Ran} F_* = \mathcal{D}_1 \oplus \mathcal{D}_2 \oplus \langle \xi_2 \rangle$, where $\mathcal{D}_1$ is invariant and $\mathcal{D}_2$ is pointwise slant with slant function $\theta_2$. In this case,*

$$\rho^{\mathcal{R}} = \frac{c+3}{4} + \frac{3(c-1)}{4r(r-1)}(r_1 + r_2 \cos^2\theta_2) - \frac{c-1}{2r}.$$

*Then*

$$\rho^{\mathcal{H}} \leq \rho^{\mathcal{R}} + \delta^{\mathcal{H}}_C(r-1) - \frac{\mathcal{C}^{\mathcal{H}}}{2(r-1)} - \frac{1}{2r}\inf_L \mathcal{C}(L) - \frac{r_2 \sin^2\theta_2}{r-1} - \frac{3\|A\|^2}{r(r-1)}. \tag{3.14}$$

*For every $\theta_2 \in [0, \pi/2]$, we also have*

$$\rho^{\mathcal{H}} \leq \rho^{\mathcal{R}} + \widehat{\delta}^{\mathcal{H}}_C(r-1) - \frac{\mathcal{C}^{\mathcal{H}}}{2(r-1)} - \frac{2r_2 \sin^2\theta_2}{r} - \frac{3\|A\|^2}{r(r-1)}. \tag{3.15}$$

*At a point where $\theta_2 > 0$ and $r_2 = \dim \mathcal{D}_2 \geq 2$, equality in* (3.14) *holds if and only if $h|_{\mathcal{D}\times\mathcal{D}} = 0$. If $0 < \theta_2 < \pi/2$, the same condition characterizes equality in* (3.15)*, with no additional rank restriction.*

**Corollary 3.5.** *Let $F : (M_1, g_1) \to (M_2(c), \varphi_2, \xi_2, \eta_2, g_2)$ be a pointwise hemi-slant Riemannian map of constant rank $r \geq 3$ into a Sasakian space form, with $\xi_2 \in \Gamma(\operatorname{Ran} F_*)$ and $\operatorname{Ran} F_* = \mathcal{D}_1 \oplus \mathcal{D}_2 \oplus \langle \xi_2 \rangle$, where $\mathcal{D}_1$ is anti-invariant and $\mathcal{D}_2$ is pointwise slant with slant function $\theta_2$. In this case,*

$$\rho^{\mathcal{R}} = \frac{c+3}{4} + \frac{3(c-1)}{4r(r-1)} r_2 \cos^2\theta_2 - \frac{c-1}{2r}.$$

*Then*

$$\rho^{\mathcal{H}} \leq \rho^{\mathcal{R}} + \delta^{\mathcal{H}}_C(r-1) - \frac{\mathcal{C}^{\mathcal{H}}}{2(r-1)} - \frac{1}{2r}\inf_L \mathcal{C}(L) - \frac{r_1 + r_2 \sin^2\theta_2}{r-1} - \frac{3\|A\|^2}{r(r-1)}. \tag{3.16}$$

*If, in addition, $r \geq 7$ and $0 < \theta_2 < \pi/2$, then*

$$\rho^{\mathcal{H}} \leq \rho^{\mathcal{R}} + \widehat{\delta}^{\mathcal{H}}_C(r-1) - \frac{\mathcal{C}^{\mathcal{H}}}{2(r-1)} - \frac{2(r_1 + r_2 \sin^2\theta_2)}{r} - \frac{2r-1}{r(r-1)}\sin^2\theta_2 - \frac{3\|A\|^2}{r(r-1)}. \tag{3.17}$$

*At a point where either $r_1 = \dim \mathcal{D}_1 \geq 2$ or $\theta_2 > 0$, equality in* (3.16) *holds if and only if $h|_{\mathcal{D}\times\mathcal{D}} = 0$. Under the additional hypotheses of the second inequality, the same condition characterizes equality in* (3.17).

**Corollary 3.6.** *Let $F : (M_1, g_1) \to (M_2(c), \varphi_2, \xi_2, \eta_2, g_2)$ be a pointwise slant Riemannian map of constant rank $r \geq 3$ into a Sasakian space form, with $\xi_2 \in \Gamma(\operatorname{Ran} F_*)$ and slant function $\theta$ on $\mathcal{D} = \operatorname{Ran} F_* \cap \ker \eta_2$. In this case,*

$$\rho^{\mathcal{R}} = \frac{c+3}{4} + \frac{3(c-1)}{4r}\cos^2\theta - \frac{c-1}{2r}.$$

*Then*

$$\rho^{\mathcal{H}} \leq \rho^{\mathcal{R}} + \delta^{\mathcal{H}}_C(r-1) - \frac{\mathcal{C}^{\mathcal{H}}}{2(r-1)} - \frac{1}{2r}\inf_L \mathcal{C}(L) - \sin^2\theta - \frac{3\|A\|^2}{r(r-1)}. \tag{3.18}$$

*If, in addition, $r \geq 7$ and $0 < \theta < \pi/2$, then*

$$\rho^{\mathcal{H}} \leq \rho^{\mathcal{R}} + \widehat{\delta}^{\mathcal{H}}_C(r-1) - \frac{\mathcal{C}^{\mathcal{H}}}{2(r-1)} - \frac{2(r-1)}{r}\sin^2\theta - \frac{2r-1}{r(r-1)}\sin^2\theta - \frac{3\|A\|^2}{r(r-1)}. \tag{3.19}$$

*At a point where $\theta > 0$, equality in* (3.18) *holds if and only if $h|_{\mathcal{D}\times\mathcal{D}} = 0$. Under the additional hypotheses of the second inequality, the same condition characterizes equality in* (3.19).

*Proof of Corollaries 3.4–3.6.* The semi-slant and hemi-slant $\delta^{\mathcal{H}}_C$ inequalities follow from Theorem A by taking, respectively, $\theta_1 = 0$ and $\theta_1 = \pi/2$. The slant inequality follows directly from the same proof, using $\|Q^{\mathcal{H}}\|^2 = (r-1)\sin^2\theta$; no decomposition of $\mathcal{D}$ into two nonzero summands is required. For the second semi-slant inequality, take a hyperplane maximizing $\mathcal{C}(L)$. The two terms inside braces in (3.9) are nonnegative, so the asserted bound follows directly, without a rank restriction and for every $\theta_2$. The other two second inequalities follow from Lemma 3.1 and the rank-dependent estimate in the proof of Theorem B for $r \geq 7$: the minimum of the two squared sines is $\sin^2\theta_2$ in the hemi-slant case and $\sin^2\theta$ in the slant case.

For the equality assertions in the first inequalities, the argument of Theorem A requires only $\operatorname{rank}(Q^{\mathcal{H}}|_{\mathcal{D}}) \geq 2$.

The additional hypotheses in the equality statements of the three corollaries ensure this rank condition. The higher-rank equality argument for the hemi-slant and slant cases is as in Theorem B. For the semi-slant $\widehat{\delta}^{\mathcal{H}}_C$ inequality, (3.9) shows that equality requires $B^\circ = 0$ and a hyperplane $L$ such that $\|h\|^2 = \|h|_{L\times L}\|^2$. If $0 < \theta_2 < \pi/2$, the restriction of $Q^{\mathcal{H}}$ to $\mathcal{D}_2$ has rank $r_2 \geq 2$. The unit normal to $L$ must therefore lie in $\mathcal{D}_1$, and the umbilical form $B$ must vanish. Conversely, if $B = 0$, a hyperplane perpendicular to a unit vector of $\mathcal{D}_1$ attains equality. No rank-seven assumption is used. □

# 4 Examples

The following examples illustrate equality in the principal inequalities, the pointwise slant specialization, and the rank-five estimate. We use the Sasakian Heisenberg space form $\mathbb{H}^{2n+1}(-3)$, represented by $\mathbb{R}^{2n+1}$ with

$$\eta_2 = dz + \sum_{j=1}^{n}(y_j\, dx_j - x_j\, dy_j), \qquad g_2 = \eta_2^2 + \sum_{j=1}^{n}(dx_j^2 + dy_j^2),$$

and the orthonormal horizontal fields

$$X_j = \frac{\partial}{\partial x_j} - y_j \frac{\partial}{\partial z}, \qquad Y_j = \frac{\partial}{\partial y_j} + x_j \frac{\partial}{\partial z}, \qquad \xi_2 = \frac{\partial}{\partial z},$$

where $\varphi_2 X_j = Y_j$ and $\varphi_2 Y_j = -X_j$. In particular, $[X_j, Y_j] = 2\xi_2$. For left-invariant horizontal fields $U, V$,

$$\nabla^{M_2}_U V = g_2(\varphi_2 U, V)\xi_2, \qquad \nabla^{M_2}_U \xi_2 = -\varphi_2 U. \tag{4.1}$$

The projection $\pi : \mathbb{H}^{2n+1}(-3) \to \mathbb{C}^n$, obtained by omitting the $z$-coordinate, is a Riemannian submersion with Reeb fibers.

**Example 4.1** (Equality in both main theorems)**.** In $\mathbb{H}^{13}(-3)$, consider the left-invariant distributions

$$\mathcal{D}_1 = \operatorname{span}\left\{X_1, \frac{Y_1 + X_2}{\sqrt{2}}\right\},$$
$$\mathcal{D}_2 = \operatorname{span}\left\{X_3, \frac{Y_3 + \sqrt{3}X_4}{2}, X_5, \frac{Y_5 + \sqrt{3}X_6}{2}\right\}.$$

Their slant angles are $\theta_1 = \pi/4$ and $\theta_2 = \pi/3$. The vector space $\mathcal{D} \oplus \langle \xi_2 \rangle$ is a Lie subalgebra of the Heisenberg algebra. Let $S$ be the corresponding connected immersed Lie subgroup, with the induced metric, and define

$$F : S \times \mathbb{R} \longrightarrow \mathbb{H}^{13}(-3), \qquad F(p, s) = p,$$

where the source has the product metric. Then $F$ is a Riemannian map of rank 7, its kernel is one-dimensional, and $A = 0$. By (4.1),

$$h(U, V) = 0 \quad (U, V \in \mathcal{D}), \qquad h(U, \xi_2) = -Q^{\mathcal{H}}U.$$

Consequently,

$$r_1 = 2, \quad r_2 = 4, \quad \|P^{\mathcal{H}}\|^2 = 2, \quad \|Q^{\mathcal{H}}\|^2 = 4,$$
$$\mathcal{C}^{\mathcal{H}} = \frac{8}{7}, \quad \rho^{\mathcal{R}} = \frac{1}{7}, \quad \rho^{\mathcal{H}} = -\frac{1}{21},$$
$$\inf_L \mathcal{C}(L) = 0, \qquad \sup_L \mathcal{C}(L) = \frac{7}{6}.$$

Thus $\delta^{\mathcal{H}}_C(6) = 4/7$ and $\widehat{\delta}^{\mathcal{H}}_C(6) = 101/84$. Substitution gives equality in both (1.1) and (1.3). In particular, these estimates are attained by a proper pointwise bi-slant map with $h \neq 0$.

**Example 4.2** (Nonintegrable horizontal distribution)**.** The construction in Example 4.1 also gives equality with $A \neq 0$. Let $\alpha = x_1\, dx_3|_S$, where $x_1, x_3$ are ambient Heisenberg coordinates restricted to $S$. Equip $S \times \mathbb{R}$ with

$$g_1 = g_S + (ds + \alpha)^2.$$

Let $F = i \circ \pi$, where $\pi : S \times \mathbb{R} \to S$ is projection and $i : S \hookrightarrow \mathbb{H}^{13}(-3)$ is the immersion above. The horizontal lifts $\widetilde{X}$ and $\widetilde{Y}$ of fields on $S$ satisfy

$$A_{\widetilde{X}}\widetilde{Y} = -\frac{1}{2}\, d\alpha(X, Y)\, \partial_s.$$

Since $d\alpha = dx_1 \wedge dx_3$ and $d\alpha(X_1, X_3) = 1$, the only nonzero $A$-components in the adapted orthonormal frame, up to skew-symmetry, are $A_{\widetilde{X}_1}\widetilde{X}_3 = -\frac{1}{2}\partial_s$ and $A_{\widetilde{X}_3}\widetilde{X}_1 = \frac{1}{2}\partial_s$. Thus $\|A\|^2 = \frac{1}{2}$. The horizontal second fundamental form of $F$ is still that of $i$, so $h|_{\mathcal{D}\times\mathcal{D}} = 0$. The Casorati invariants and $\rho^{\mathcal{R}}$ therefore have the same values as above, while (3.5) gives

$$\rho^{\mathcal{H}} = -\frac{1}{21} - \frac{3(1/2)}{7 \cdot 6} = -\frac{1}{12}.$$

Both right-hand sides in (1.1) and (1.3) are their preceding values minus 1/28, hence equal $-1/12$. Thus both inequalities are attained with $A \neq 0$.

**Example 4.3** (Equality in the pointwise slant corollary)**.** In the Sasakian space form $\mathbb{H}^{13}(-3)$, consider the orthonormal horizontal vector fields

$$E_{2j-1} = X_{2j-1}, \qquad E_{2j} = \frac{Y_{2j-1} + X_{2j}}{\sqrt{2}}, \qquad j = 1, 2, 3.$$

Put $\mathcal{D} = \operatorname{span}\{E_1, \dots, E_6\}$. Since

$$P^{\mathcal{H}} E_{2j-1} = \frac{1}{\sqrt{2}} E_{2j}, \qquad P^{\mathcal{H}} E_{2j} = -\frac{1}{\sqrt{2}} E_{2j-1},$$

the distribution $\mathcal{D}$ is slant with constant angle $\theta = \pi/4$. The vector space $\mathcal{D} \oplus \langle \xi_2 \rangle$ is a Lie subalgebra of the Heisenberg algebra. Let $S$ be the corresponding connected immersed Lie subgroup with the induced metric, and define

$$F : S \times \mathbb{R} \longrightarrow \mathbb{H}^{13}(-3), \qquad F(p, s) = p,$$

where the source carries the product metric. Then $F$ is a Riemannian map of rank 7 with nontrivial kernel and $A = 0$. By (4.1),

$$h(U, V) = 0, \qquad U, V \in \mathcal{D},$$

whereas

$$h(U, \xi_2) = -Q^{\mathcal{H}} U, \qquad h(\xi_2, \xi_2) = 0.$$

Consequently,

$$\|P^{\mathcal{H}}\|^2 = 3, \quad \|Q^{\mathcal{H}}\|^2 = 3, \quad \|h\|^2 = 6, \quad \operatorname{tr} h = 0,$$
$$\mathcal{C}^{\mathcal{H}} = \frac{6}{7}, \quad \rho^{\mathcal{R}} = \frac{1}{14}, \quad \rho^{\mathcal{H}} = -\frac{1}{14}.$$

The hyperplane $L = \mathcal{D}$ satisfies $\mathcal{C}(L) = 0$. To determine the supremum, let an arbitrary hyperplane have unit normal $\sqrt{t}\,\xi_2 + \sqrt{1-t}\,u$, where $u \in \mathcal{D}$ is a unit vector and $0 \le t \le 1$. Using the preceding expressions for $h$, we obtain

$$\|h\|^2 - \|h|_{L\times L}\|^2 = 1 + 3t + 2t^2 \geq 1.$$

Equality holds for $t = 0$. Hence

$$\inf_L \mathcal{C}(L) = 0, \qquad \sup_L \mathcal{C}(L) = \frac{5}{6}.$$

It follows that

$$\delta^{\mathcal{H}}_C(6) = \frac{3}{7}, \qquad \widehat{\delta}^{\mathcal{H}}_C(6) = \frac{79}{84}.$$

Substituting these values into (3.18) gives

$$-\frac{1}{14} = \frac{1}{14} + \frac{7}{8}\cdot\frac{3}{7} - \frac{1}{56} - \frac{1}{2}.$$

Similarly, (3.19) gives

$$-\frac{1}{14} = \frac{1}{14} + \frac{79}{84} - \frac{1}{14} - \frac{6}{7} - \frac{13}{84}.$$

Thus equality is attained in both inequalities of the pointwise slant corollary, although $h \neq 0$. In contrast, the corresponding general $\delta_C$ estimate in [25] gives $\rho^{\mathcal{H}} \le \rho^{\mathcal{R}} + \delta^{\mathcal{H}}_C(6) = \frac{1}{2}$, whereas the refined bound and the actual curvature both equal $-\frac{1}{14}$.

**Example 4.4** (The rank-five case)**.** We first observe that the first alternative of Proposition 3.2 is attained with $h|_{\mathcal{D}\times\mathcal{D}} = 0$. Indeed, in $\mathbb{H}^9(-3)$ choose the left-invariant tangent space generated by $\xi_2$ and

$$X_1, \quad \frac{Y_1 + X_2}{\sqrt{2}}, \quad X_3, \quad \frac{Y_3 + 3X_4}{\sqrt{10}}.$$

The corresponding Reeb-invariant Lie subgroup, followed by the product construction used above, gives a Riemannian map of rank 5 with $\sin^2\theta_1 = 1/2$ and $\sin^2\theta_2 = 9/10$. Here

$$\rho^{\mathcal{R}} = \frac{11}{50}, \qquad \rho^{\mathcal{H}} = -\frac{3}{50}, \qquad \mathcal{C}^{\mathcal{H}} = \frac{28}{25}, \qquad \widehat{\delta}^{\mathcal{H}}_C(4) = \frac{241}{200},$$

and equality holds in the first alternative of Proposition 3.2.

To realize the second alternative, identify the horizontal space of $\mathbb{H}^9(-3)$ at the identity with $\mathbb{C}^4$. Choose the orthonormal vectors

$$E_1 = X_1, \qquad E_2 = \frac{Y_1 + X_2}{\sqrt{2}}, \qquad E_3 = X_3, \qquad E_4 = \sqrt{\frac{2}{5}}\,Y_3 + \sqrt{\frac{3}{5}}\,X_4,$$

regarded as constant vectors in $\mathbb{C}^4$ when defining the graph below. Their two orthogonal slant planes have

$$\sin^2\theta_1 = \frac{1}{2}, \qquad \sin^2\theta_2 = \frac{3}{5} < \frac{13}{8}\sin^2\theta_1.$$

Let $J$ denote the standard complex structure of $\mathbb{C}^4$, identified with the restriction of $\varphi_2$ to the Heisenberg horizontal space. Let $N_i$ be the normal component of $JE_i$ relative to

span$\{E_1, E_2, E_3, E_4\} \subset \mathbb{C}^4$. The $N_i$ are mutually orthogonal, with $\|N_1\|^2 = \|N_2\|^2 = 1/2$ and $\|N_3\|^2 = \|N_4\|^2 = 3/5$. Set

$$t = \frac{17\sqrt{5} - 34}{9\sqrt{5} + 34}, \qquad a = \frac{2t^{3/2}}{(1+t)\sqrt{1-t}}, \qquad b = \frac{9\sqrt{t(1-t)}}{17 - 9t}.$$

Consider, near the origin of $\mathbb{R}^4$, the immersion

$$f(x_1, x_2, x_3, x_4) = \sum_{i=1}^{4} x_i E_i + \frac{a}{2} \sum_{i=1}^{4} x_i^2 N_1 + b x_1 \sum_{j=2}^{4} x_j N_j. \tag{4.2}$$

Let $S$ be the inverse image of $f(\mathbb{R}^4)$ under the Reeb projection, restricted to a sufficiently small neighborhood, and let $F : S \times \mathbb{R} \to \mathbb{H}^9(-3)$ be the projection onto $S$ followed by the inclusion. The two distinct proper slant planes at the origin extend locally to smooth pointwise slant distributions, namely the two rank-two eigendistributions of $-(P^{\mathcal{H}})^2$ on the range orthogonal to $\xi_2$. Thus $F$ is locally a proper pointwise bi-slant Riemannian map.

At the origin, with $B = h|_{\mathcal{D}\times\mathcal{D}}$, the graph (4.2) gives

$$B(E_i, E_i) = aN_1 \quad (1 \le i \le 4), \qquad B(E_1, E_j) = bN_j \quad (2 \le j \le 4),$$

and all remaining off-diagonal components vanish. Moreover, $h(E_i, \xi_2) = -N_i$, $h(\xi_2, \xi_2) = 0$ and $A = 0$. In particular,

$$\|h\|^2 = \frac{22}{5} + 2a^2 + \frac{17}{5}b^2, \qquad \|\operatorname{tr} h\|^2 = 8a^2.$$

Let $L$ be the hyperplane orthogonal to $\sqrt{1-t}\,E_1 + \sqrt{t}\,\xi_2$. Direct substitution yields

$$\begin{aligned} &\frac{9}{2}\Big(\|h\|^2 - \|h|_{L\times L}\|^2\Big) + 4\|B^\circ\|^2 \\ &= \frac{36}{13}\left[\sqrt{17 \cdot \frac{1}{2}\left(\frac{1}{2} + 2 \cdot \frac{3}{5}\right)} - 2\left(\frac{1}{2} + \frac{3}{5}\right)\right]. \end{aligned} \tag{4.3}$$

The lower bound in the proof of Proposition 3.2 holds for every hyperplane; hence $L$ realizes the supremum defining $\widehat{\delta}_C^{\mathcal{H}}(4)$. Equations (3.9) and (4.3) therefore show that equality is attained in the second alternative of Proposition 3.2 as well.

# 5 Conclusion

We have established sharp Casorati inequalities for pointwise bi-slant Riemannian maps into Sasakian space forms with the characteristic vector field tangent to the range. The Sasakian identities governing the mixed components of the second fundamental form yield refined bounds for the normalized horizontal scalar curvature while retaining the contribution of O'Neill's tensor. We have characterized the equality cases of the principal inequalities and obtained a separate sharp estimate in rank five. The corresponding inequalities for pointwise semi-slant, hemi-slant, and slant Riemannian maps have also been derived. Examples in Sasakian Heisenberg space forms demonstrate that equality can occur even when the second fundamental form is nonzero and the horizontal distribution is nonintegrable.